\documentclass[11pt]{amsart}
\usepackage{amssymb}
\usepackage{amsmath}
\usepackage[active]{srcltx}
\usepackage{t1enc}
\usepackage[latin2]{inputenc}
\usepackage{verbatim}
\usepackage{amsmath,amsfonts,amssymb,amsthm}
\usepackage[mathcal]{eucal}
\usepackage{enumerate}
\usepackage[centertags]{amsmath}
\usepackage{graphics}

\newtheorem{theorem}{Theorem}

\newtheorem{lemma}{Lemma}

\begin{document}
\author{I. BLAHOTA AND G. TEPHNADZE}
\title[On the $\left(C,\alpha \right) $-means\dots]{On the $\left(C,\alpha
\right) $-means with respect to the Walsh system}
\address{I. Blahota, Institute of Mathematics and Computer Sciences, College
of Ny\'\i regyh\'aza, P.O. Box 166, Ny\'\i regyh\'aza, H-4400, Hungary}
\email{blahota@nyf.hu}
\address{G. Tephnadze, Department of Mathematics, Faculty of Exact and
Natural Sciences, Tbilisi State University, Chavchavadze str. 1, Tbilisi
0128, Georgia}
\email{giorgitephnadze@gmail.com}
\thanks{The first author was supported by project TÁ%
MOP-4.2.2.A-11/1/KONV-2012-0051.\\
The second author was supported by Shota Rustaveli National Science
Foundation grant no. 52/54 (Bounded operators on the martingale Hardy
spaces).}
\date{}
\maketitle

\begin{abstract}
As main result we prove strong convergence theorems for Ces\'aro means $%
\left(C,\alpha \right) $ on the Hardy spaces $H_{1/\left(1+\alpha \right) } $%
, where $0<\alpha <1.$
\end{abstract}

\date{}



\section{INTRODUCTION}

It is well-known that Walsh-Paley system does not form a basis in the space $%
L_{1}\left(G\right).$ Moreover, there is a function in the dyadic Hardy
space $H_{1}\left(G\right) ,$ such that the partial sums of $F$ are not
bounded in $L_{1}$-norm. However, in Simon \cite{Si3} the following
estimation was obtained for all $F\in H_{1}\left(G\right):$ 
\begin{equation*}
\frac{1}{\log n}\overset{n}{\underset{k=1}{\sum }}\frac{\left\Vert
S_{k}F\right\Vert _{1}}{k}\leq c\left\Vert F\right\Vert _{H_{1}},
\end{equation*}
where $S_{k}F$ denotes the $k$-th partial sum of the Walsh-Fourier series of 
$F.$ (For the trigonometric analogue see in Smith \cite{sm}, for the
Vilenkin system in G\'at \cite{gat1}). Simon \cite{si1} (see also \cite{We})
proved that there is an absolute constant $c_{p},$ depending only on $p,$
such that 
\begin{equation}
\frac{1}{\log ^{\left[ p\right] }n}\overset{n}{\underset{k=1}{\sum }}\frac{
\left\Vert S_{k}F\right\Vert _{p}^{p}}{k^{2-p}}\leq c_{p}\left\Vert
F\right\Vert _{H_{p}}^{p},\text{ \ \ }\left(0<p\leq 1\right) ,  \label{1cc}
\end{equation}
for all $F\in H_{p}$ and $n\in \mathbb{N}$, where $\left[ p\right] $ denotes
integer part of $p.$

The second author \cite{tep4} proved that sequence $\left\{
1/k^{2-p}\right\} _{k=1}^{\infty }$ $\left(0<p<1\right) $ \ in (\ref{1cc})
is given exactly.

Weisz \cite{We3} considered the norm convergence of Fej\'er means of
Walsh-Fourier series and proved that the following is true:

\textbf{Theorem W1. }Let $F\in H_{p}.$ Then 
\begin{equation}
\left\Vert \sigma _{k}F\right\Vert _{H_{p}}\leq c_{p}\left\Vert F\right\Vert
_{H_{p}},\text{ \ \ \ \ \ \ }\left(1/2<p<\infty \right)  \label{1ac}
\end{equation}

This theorem implies that 
\begin{equation*}
\frac{1}{n^{2p-1}}\overset{n}{\underset{k=1}{\sum }}\frac{\left\Vert \sigma
_{k}F\right\Vert _{H_{p}}^{p}}{k^{2-2p}}\leq c_{p}\left\Vert F\right\Vert
_{H_{p}}^{p},\text{ \ \ \ }\left(1/2<p<\infty \right) .
\end{equation*}

If Theorem W1 held for $0<p\leq 1/2,$ then we would have 
\begin{equation}
\frac{1}{\log ^{\left[ 1/2+p\right] }n}\overset{n}{\underset{k=1}{\sum }} 
\frac{\left\Vert \sigma _{k}F\right\Vert _{H_{p}}^{p}}{k^{2-2p}}\leq
c_{p}\left\Vert F\right\Vert _{H_{p}}^{p},\text{ \ \ \ }\left(0<p\leq
1/2\right) ,  \label{2cc}
\end{equation}

\bigskip but the second author \cite{tep1} proved that the assumption $p>1/2$
is essential. In particular, he proved that there exists a martingale $F\in
H_{p}$ $\left(0<p\leq 1/2\right) ,$ such that 
\begin{equation*}
\sup_{n}\left\Vert \sigma _{n}F\right\Vert _{p}=+\infty .
\end{equation*}

However, the second author \cite{tep5} prove that (\ref{2cc}) holds, though
( \ref{1ac}) is not true for $0<p\leq 1/2.$

The weak (1,1)-type inequality for the maximal operator of Fej\'er means 
\begin{equation*}
\mu \left(\sigma ^{\ast }f>\lambda \right) \leq \frac{c}{\lambda }
\left\Vert f\right\Vert _{1},\text{ \qquad }\left(\lambda >0\right)
\end{equation*}
can be found in Schipp \cite{Sc} (see also \cite{PS}). Fujji \cite{Fu} and
Simon \cite{Si2} verified that $\sigma ^{\ast }$ is bounded from $H_{1}$ to $%
L_{1}$. Weisz \cite{We2} generalized this result and proved the boundedness
of $\sigma ^{\ast }$ from the martingale space $H_{p}$ to the space $L_{p}$
for $p>1/2$. Simon \cite{Si1} gave a counterexample, which shows that
boundedness does not hold for $0<p<1/2.$ The counterexample for $p=1/2$ due
to Goginava \cite{GoAMH}, (see also \cite{BGG2} and \cite{tep1}). Weisz \cite%
{we4} proved that $\sigma ^{\ast }$ is bounded from the Hardy space $H_{1/2}$
to the space $L_{1/2,\infty }$.

The second author \cite{tep2, tep3} proved that the following is true:

\textbf{Theorem T1.} The maximal operators $\widetilde{\sigma }_{p}^{\ast }$
defined by 
\begin{equation}
\widetilde{\sigma }_{p}^{\ast }:=\sup_{n\in \mathbb{N}}\frac{\left\vert
\sigma _{n}\right\vert }{n^{1/p-2}\log ^{2\left[ 1/2+p\right] }n},\text{ \ }
\left(0<p\leq 1/2,\text{ }n=2,3,...\right)  \label{1}
\end{equation}
where $\left[ 1/2+p\right] $ denotes integer part of $1/2+p,$ is bounded
from the Hardy space $H_{p}$ to the space $L_{p}.$ Moreover, there was also
shown that sequence $\left\{ n^{1/p-2}\log ^{2\left[ 1/2+p\right]
}n:n=2,3,...\right\} $ in (\ref{1}) can not be improved.

The maximal operator $\sigma ^{\alpha ,\ast }$ $\left(0<\alpha <1\right) $
of the Ces\'aro means means of Walsh-Paley system was investigated by Weisz 
\cite{we6}. In his paper Weisz proved that $\sigma ^{\alpha ,\ast }$ is
bounded from the martingale space $H_{p}$ to the space $L_{p}$ for $%
p>1/\left(1+\alpha \right) .$ Goginava \cite{gog4} gave counterexample,
which shows that boundedness does not hold for $0<p\leq 1/\left(1+\alpha
\right) .$ Recently, Weisz and Simon \cite{sw} show that in case $p=1/\left(
1+\alpha \right) $ the maximal operator $\sigma ^{\alpha ,\ast }$ is bounded
from the Hardy space $H_{1/\left(1+\alpha \right) }$ to the space $%
L_{1/\left(1+\alpha \right) ,\infty }$.

In \cite{gog8} Goginava investigated the behaviour of Ces\'aro means of
Walsh-Fourier series in detail. For some approximation properties of the two
dimensional case see paper of Nagy \cite{na}.

The main aim of this paper is to generalize Theorem T1 and estimation (\ref%
{2cc}) for Ces\'aro means, when $p=1/\left(1+\alpha \right) .$

\section{Definitions and Notations}

Let $\mathbb{N}_{+}$ denote the set of the positive integers, $\mathbb{N}:= 
\mathbb{N}_{+}\cup \{0\}.$ Denote by $Z_{2}$ the discrete cyclic group of
order 2, that is $Z_{2}:=\{0,1\},$ where the group operation is the modulo 2
addition and every subset is open. The Haar measure on $Z_{2}$ is given so
that the measure of a singleton is 1/2.

Define the group $G$ as the complete direct product of the group $Z_{2}$
with the product of the discrete topologies of $Z_{2}`$s.

The elements of $G$ are represented by sequences 
\begin{equation*}
x:=(x_{0},x_{1},\dots ,x_{j},\dots )\qquad \left(x_{k}\in \{0,1\}\right) .
\end{equation*}

It is easy to give a base for the neighbourhood of $G$ 
\begin{equation*}
I_{0}\left(x\right) :=G,
\end{equation*}
\begin{equation*}
I_{n}(x):=\{y\in G\mid y_{0}=x_{0},\dots ,y_{n-1}=x_{n-1}\}\text{ }(x\in G, 
\text{ }n\in \mathbb{N}).
\end{equation*}
Denote $I_{n}:=I_{n}\left(0\right) $ for $n\in \mathbb{N}$ and $\overline{
I_{n}}:=G$ $\backslash $ $I_{n}$. Let

\begin{equation*}
e_{n}:=\left(0,\dots,0,x_{n}=1,0,\dots\right) \in G\qquad \left(n\in \mathbb{%
N}\right)
\end{equation*}

It is evident 
\begin{equation}
\overline{I_{M}}=\left(\overset{M-2}{\underset{k=0}{\bigcup }}\overset{M-1}{ 
\underset{l=k+1}{\bigcup }}I_{l+1}\left(e_{k}+e_{l}\right) \right) \bigcup
\left(\underset{k=0}{\bigcup\limits^{M-1}}I_{M}\left(e_{k}\right) \right) .
\label{2}
\end{equation}

The norm (or quasi-norm) of the space $L_{p}(G)$ is defined by \qquad

\begin{equation*}
\left\| f\right\| _{p}:=\left(\int_{G}\left| f(x)\right| ^{p}d\mu (x)\right)
^{1/p},\qquad \left(0<p<\infty \right).
\end{equation*}

The space $L_{p,\infty }\left(G\right) $ consists of all measurable
functions $f$ for which

\begin{equation*}
\left\Vert f\right\Vert _{L_{p,\infty }(G)}:=\underset{\lambda >0}{\sup }
\lambda \mu \left(f>\lambda \right) ^{1/p}<+\infty .
\end{equation*}

If $n\in \mathbb{N}$ then for every $n$ can be uniquely expressed as $%
n=\sum_{k=0}^{\infty }n_{j}2^{j}$ where $n_{j}\in Z_{2}$ $~(j\in \mathbb{N})$
and only a finite number of $n_{j}$'s differs from zero. Let $\left\vert
n\right\vert :=\max $ $\{j\in \mathbb{N},$ $n_{j}\neq 0\},$ that is $%
2^{\left\vert n\right\vert }\leq n\leq 2^{\left\vert n\right\vert +1}.$

Next, we introduce on $G$ an orthonormal system which is called the Walsh
system. At first define the Rademacher functions as 
\begin{equation*}
r_{k}\left(x\right) :=\left(-1\right) ^{x_{k}}\text{\qquad } \left(x\in G,%
\text{ }k\in \mathbb{N}\right) .
\end{equation*}

Now define the Walsh system $w:=(w_{n}:n\in \mathbb{N})$ on $G$ as: 
\begin{equation*}
w_{n}(x):=\prod^{\infty}_{k=0}r_{k}^{n_{k}}\left(x\right) =r_{\left\vert
n\right\vert }\left(x\right) \left(-1\right) ^{\underset{k=0}{\overset{
\left\vert n\right\vert -1}{\sum }}n_{k}x_{k}}\text{\qquad }\left(n\in 
\mathbb{N}\right).
\end{equation*}

The Walsh system is orthonormal and complete in $L_{2}\left( G\right) $ (see 
\cite{sws}).

If $f\in L_{1}\left(G\right) $ we can establish the Fourier coefficients,
the partial sums of the Fourier series, the Fej\'er means, the Dirichlet and
Fej\'er kernels in the usual manner: 
\begin{eqnarray*}
\widehat{f}\left(k\right) &:&=\int_{G}fw_{k}d\mu \,\,\,\,\qquad \left(k\in 
\mathbb{N}\right) , \\
S_{n}f &:&=\sum_{k=0}^{n-1}\widehat{f}\left(k\right) w_{k}\text{ \qquad }
\left(n\in \mathbb{N}_{+},\ S_{0}f:=0\right) , \\
\sigma _{n}f &:&=\frac{1}{n}\sum_{k=1}^{n}S_{k}f\text{\qquad }\left(n\in 
\mathbb{N}_{+}\right) , \\
D_{n} &:&=\sum_{k=0}^{n-1}w_{k\text{ }}\,\,\qquad \,\left(n\in \mathbb{N}
_{+}\right) , \\
K_{n} &:&=\frac{1}{n}\overset{n}{\underset{k=1}{\sum }}D_{k}\text{ \qquad
\thinspace }\left(n\in \mathbb{N}_{+}\right) ,
\end{eqnarray*}
respectively.

Recall that 
\begin{equation}
D_{2^{n}}\left(x\right) =\left\{ 
\begin{array}{ll}
2^{n}, & \text{if\thinspace \thinspace \thinspace }x\in I_{n}, \\ 
0, & \text{if}\,\,x\notin I_{n}.%
\end{array}
\right.  \label{1dn}
\end{equation}

For the $2^{n}$-th Fej\'er kernel we have the following equality (see \cite%
{gat}): 
\begin{equation}
K_{2^{n}}\left(x\right) =\left\{ 
\begin{array}{ll}
2^{t-1}, & \text{if } x\in I_{n}(e_{t}), \\ 
\left(2^{A}+1\right)/2, & \text{if } x\in I_{n}, \\ 
0, & \text{otherwise.}%
\end{array}
\right.  \label{5a}
\end{equation}
for $n>t,$ $t,n\in \mathbb{N},$ $x\in I_{t}\backslash $ $I_{t+1}$.

The Ces\'aro means, ($\left(C,\alpha \right) $ means) and it`s kernel with
respect to the Walsh-Fourier series are defined as

\begin{equation*}
\sigma _{n}^{\alpha }f:=\frac{1}{A_{n}^{\alpha }}\overset{n}{\underset{k=1}{
\sum }}A_{n-k}^{\alpha -1}S_{k}f,\text{ \ }K_{n}^{\alpha }f:=\frac{1}{
A_{n}^{\alpha }}\overset{n}{\underset{k=1}{\sum }}A_{n-k}^{\alpha -1}D_{k}f,
\end{equation*}
respectively, where

\begin{equation}
A_{0}^{\alpha }:=0,\text{ \qquad }A_{n}^{\alpha }:=\frac{\left(\alpha
+1\right) \dots\left(\alpha +n\right) }{n!},\ \alpha \neq -1,-2,\dots
\label{1d}
\end{equation}

It is well known that

\begin{equation*}
A_{n}^{\alpha }=\overset{n}{\underset{k=0}{\sum }}A_{n-k}^{\alpha -1},\text{
\ \ \ \ \ \ }A_{n}^{\alpha }-A_{n-1}^{\alpha }=A_{n}^{\alpha -1},\text{ \ \
\ \ \ \ }A_{n}^{\alpha }\backsim n^{\alpha },
\end{equation*}

and 
\begin{equation}
\sup_{n}\int_{G}\left\vert K_{n}^{\alpha }\left(x\right) \right\vert d\mu
\left(x\right) \leq c<\infty .  \label{4}
\end{equation}

The $\sigma $-algebra is generated by the intervals $\left\{ I_{n}\left(
x\right) :x\in G\right\} $ will be denoted by $\digamma _{n}\left(n\in 
\mathbb{N}\right) .$ The conditional expectation operators relative to $%
\digamma _{n}\left(n\in \mathbb{N}\right) $ are denoted by $E_{n}.$

A sequence $F=\left(F_{n},\text{ }n\in \mathbb{N}\right) $ of functions $%
F_{n}\in L_{1}\left(G\right) $ is said to be a dyadic martingale if (for
details see e.g. \cite{We1})

$\left(i\right) $ $F_{n}$ is $\digamma _{n}$ measurable for all $n\in 
\mathbb{N},$

$\left(ii\right) $ $E_{n}F_{m}=F_{n}$ for all $n\leq m.$

The maximal function of a martingale $F$ is defined by

\begin{equation*}
F^{\ast }=\sup_{n\in \mathbb{N}}\left\vert F_{n}\right\vert.
\end{equation*}

In case of $f\in L_{1}\left(G\right),$ the maximal functions are also be
given by

\begin{equation*}
f^{\ast }\left(x\right) =\sup\limits_{n\in \mathbb{N}}\frac{1}{\mu \left(
I_{n}\left(x\right) \right) }\left\vert \int\limits_{I_{n}\left(x\right)
}f\left(u\right) d\mu \left(u\right) \right\vert.
\end{equation*}

For $0<p<\infty $ \ the Hardy martingale spaces $H_{p}$ $\left(G\right) $
consist of all martingales for which

\begin{equation*}
\left\Vert F\right\Vert _{H_{p}}:=\left\Vert F^{\ast }\right\Vert
_{p}<\infty.
\end{equation*}

A bounded measurable function $a$ is a $p$-atom, if there exists a dyadic
interval $I$, such that \qquad 
\begin{equation*}
\int_{I}ad\mu =0,\text{ \ \ \ \ \ }\left\Vert a\right\Vert _{\infty }\leq
\mu \left(I\right) ^{-1/p},\text{ \ \ \ supp}\left(a\right) \subset I.
\end{equation*}

The dyadic Hardy martingale spaces $H_{p}$ $\left(G\right) $ for $0<p\leq 1$
have an atomic characterization. Namely, the following theorem is true (see 
\cite{We5}):

\textbf{Theorem W}: A martingale $F=\left(F_{n},\text{ }n\in \mathbb{N}
\right) $ is in $H_{p}\left(0<p\leq 1\right) $ if and only if there exists a
sequence $\left(a_{k},\text{ }k\in \mathbb{N}\right) $ of $p$-atoms and a
sequence $\left(\mu _{k},\text{ }k\in \mathbb{N}\right) $ of a real numbers
such that for every $n\in \mathbb{N}$

\begin{equation}
\qquad \sum_{k=0}^{\infty }\mu _{k}S_{2^{n}}a_{k}=F_{n}  \label{2A}
\end{equation}
and

\begin{equation*}
\qquad \sum_{k=0}^{\infty }\left\vert \mu _{k}\right\vert ^{p}<\infty ,
\end{equation*}
Moreover, $\left\Vert F\right\Vert _{H_{p}}\backsim \inf \left(
\sum_{k=0}^{\infty }\left\vert \mu _{k}\right\vert ^{p}\right) ^{1/p},$
where the infimum is taken over all decompositions of $F$ of the form (\ref%
{2A}).

It is easy to check that for every martingales $F=\left(F_{n},n\in \mathbb{N 
}\right) $ and every $k\in \mathbb{N}$ the limit

\begin{equation}
\widehat{F}\left(k\right) :=\lim_{n\rightarrow \infty }\int_{G}F_{n}\left(
x\right) w_{k}\left(x\right) d\mu \left(x\right)  \label{3a}
\end{equation}
exists and it is called the $k$-th Walsh-Fourier coefficients of $F.$

If $F:=$ $\left(E_{n}f:n\in \mathbb{N}\right) $ is a regular martingale
generated by $f\in L_{1}\left(G\right),$ then \qquad \qquad \qquad \qquad

\begin{equation*}
\widehat{F}\left(k\right) =\int_{G}f\left(x\right) w_{k}\left(x\right) d\mu
\left(x\right) =:\widehat{f}\left(k\right) ,\text{ \qquad }k\in \mathbb{N}.
\end{equation*}
For $0<\alpha <1$ let consider maximal operators

\begin{equation*}
\sigma ^{\alpha ,\ast }F:=\sup_{n\in \mathbb{N}}\left\vert \sigma
_{n}^{\alpha }F\right\vert ,\text{ \ \ \ \ }\overset{\sim }{\sigma }^{\alpha
,\ast }F:=\sup_{n\in \mathbb{N}}\frac{\left\vert \sigma _{n}^{\alpha
}F\right\vert }{\log ^{1+\alpha }n},\text{ }
\end{equation*}
For the martingale 
\begin{equation*}
F=\overset{\infty }{\underset{n=0}{\sum }}\left(F_{n}-F_{n-1}\right)
\end{equation*}
the conjugate transforms are defined as 
\begin{equation*}
\widetilde{F^{\left(t\right) }}=\overset{\infty }{\underset{n=0}{\sum }}
r_{n}\left(t\right) \left(F_{n}-F_{n-1}\right) ,
\end{equation*}
where $t\in G$ is fixed. We note that $\widetilde{F^{\left(0\right) }}=F.$
As it is well known (see \cite{We1}) 
\begin{equation}
\left\Vert \widetilde{F^{\left(t\right) }}\right\Vert _{H_{p}}=\left\Vert
F\right\Vert _{H_{p}},\text{ \ \ }\left\Vert F\right\Vert _{H_{p}}^{p}\sim
\int_{G}\left\Vert \widetilde{F^{\left(t\right) }}\right\Vert _{p}^{p}dt.
\label{5.1}
\end{equation}

\section{Formulation of Main Results}

\begin{theorem}
a) Let $0<\alpha <1$ and $f\in H_{1/\left( 1+\alpha \right) }.$ Then there
exists absolute constant $c_{\alpha },$ depending only on $\alpha ,$ such
that 
\begin{equation*}
\left\Vert \overset{\sim }{\sigma }^{\alpha ,\ast }F\right\Vert
_{H_{1/\left( 1+\alpha \right) }}\leq c_{\alpha }\left\Vert F\right\Vert
_{H_{1/\left( 1+\alpha \right) }}.
\end{equation*}%
b) Let $0<\alpha <1$ and $\varphi :\mathbb{N}_{+}\rightarrow \lbrack
1,\infty )$ be a non-decreasing function satisfying the condition 
\begin{equation}
\overline{\lim_{n\rightarrow \infty }}\frac{\log ^{1+\alpha }n}{\varphi
\left( n\right) }=+\infty ,  \label{6}
\end{equation}%
\textit{then there exists a martingale} $f\in H_{1/(1+\alpha )}\left(
G\right) ,$ \textit{such that} 
\begin{equation*}
\sup_{n\in \mathbb{N}}\left\Vert \frac{\sigma _{n}^{\alpha }f}{\varphi
\left( n\right) }\right\Vert _{1/\left( 1+\alpha \right) }=\infty .
\end{equation*}
\end{theorem}

\begin{theorem}
Let $0<\alpha <1$ and $f\in H_{1/\left( 1+\alpha \right) }..$ Then there
exists an absolute constant $c_{\alpha },$ depending only on $\alpha ,$ such
that 
\begin{equation*}
\frac{1}{\log n}\overset{n}{\underset{m=1}{\sum }}\frac{\left\Vert \sigma
_{m}^{\alpha }F\right\Vert _{H_{1/\left( 1+\alpha \right) }}^{1/\left(
1+\alpha \right) }}{m}\leq c_{\alpha }\left\Vert F\right\Vert _{H_{1/\left(
1+\alpha \right) }}^{1/\left( 1+\alpha \right) }.
\end{equation*}
\end{theorem}

\section{AUXILIARY PROPOSITIONS}

\begin{lemma}
\cite{We1} Suppose that an operator $T$ is $\sigma $-sub-linear and for some 
$0<p\leq 1$
\end{lemma}

\begin{equation*}
\int\limits_{\overset{-}{I}}\left\vert Ta\right\vert ^{p}d\mu \leq
c_{p}<\infty ,
\end{equation*}
\textit{for every} $p$\textit{-atom} $a$, \textit{where} $I$ \textit{denotes
the support of the atom. If} $T$ \textit{is bounded from} $L_{\infty \text{ }
}$ \textit{to } $L_{\infty },$ \textit{then} 
\begin{equation*}
\left\Vert Tf\right\Vert _{L_{p}\left(G\right) }\leq c_{p}\left\Vert
f\right\Vert _{H_{p}\left(G\right) }.
\end{equation*}

\begin{lemma}
\cite{gg} Let $0<\alpha <1.$ Then 
\begin{equation*}
\left\vert K_{n}^{\alpha }\right\vert \leq \frac{c_{\alpha }}{
A_{n-1}^{\alpha }}\left\{ \sum_{j=0}^{\left\vert n\right\vert }2^{j\alpha
}K_{2^{j}}\right\} ,
\end{equation*}
where $K_{n}$ and $K_{n}^{\alpha }$ are kernels of Fej\'er and Ces\'aro
means, respectively.
\end{lemma}

\begin{lemma}
Let $0<\alpha <1$ and $n>2^{M}.$ Then 
\begin{equation*}
\int_{I_{M}}\left\vert K_{n}^{\alpha }\left( x+t\right) \right\vert d\mu
\left( t\right) \leq \frac{c_{\alpha }2^{\alpha l+k}}{n^{\alpha }2^{M}},%
\text{ \ for }x\in I_{l+1}\left( e_{k}+e_{l}\right) ,
\end{equation*}%
$\text{ }(k=0,\dots ,M-2,$ $l=k+1,\dots ,M-1)$ and 
\begin{equation*}
\int_{I_{M}}\left\vert K_{n}^{\alpha }\left( x+t\right) \right\vert d\mu
\left( t\right) \leq \frac{c_{\alpha }2^{k}}{2^{M}},\text{ \ \textit{for} \ }%
x\in I_{M}\left( e_{k}\right) ,\text{ }(k=0,\dots ,M-1).
\end{equation*}
\end{lemma}

\textbf{Proof.} Let $x\in I_{l+1}\left(e_{k}+e_{l}\right) .$ Then applying (%
\ref{5a}) we have 
\begin{equation*}
K_{2^{A}}\left(x\right) =0,\text{ \qquad when \ \ }A>l.
\end{equation*}

Suppose that $k<A\leq l.$ Using (\ref{5a}) we get 
\begin{equation}
\left\vert K_{2^{A}}\left(x\right) \right\vert \leq c2^{k}.  \notag
\end{equation}

Let $A\leq k<l.$ Then

\begin{equation}
\left\vert K_{2^{A}}\left(x\right) \right\vert =\left\vert K_{2^{A}}\left(
0\right) \right\vert =\frac{2^{A}+1}{2}\leq c2^{k}.  \label{6b}
\end{equation}

If we apply Lemma 2 we conclude that

\begin{equation}
A_{n}^{\alpha }\left\vert K_{n}^{\alpha }\left(x\right) \right\vert \leq
c_{\alpha }\overset{l-1}{\underset{A=0}{\sum }}2^{\alpha A}\left\vert
K_{2^{A}}\left(x\right) \right\vert \leq c_{\alpha }\overset{l-1}{ \underset{%
A=0}{\sum }}2^{\alpha A+k}\leq c2^{\alpha l+k}.  \label{7}
\end{equation}

Let $x\in I_{l+1}\left( e_{k}+e_{l}\right) ,$ for some $0\leq k<l\leq M-1.$
Since $x+t\in $ $I_{l+1}\left( e_{k}+e_{l}\right) ,$ for $t\in $ $I_{M}$ and 
$n\geq 2^{M}$ from (\ref{7}) we obtain 
\begin{equation}
\int_{I_{M}}\left\vert K_{n}^{\alpha }\left( x+t\right) \right\vert d\mu
\left( t\right) \leq \frac{c_{\alpha }2^{\alpha l+k}}{n^{\alpha }2^{M}}.
\label{8}
\end{equation}

Let $x\in I_{M}\left(e_{k}\right) ,$ $k=0,\dots ,M-1,$ then applying Lemma 2
and (\ref{5a}) we have

\begin{equation}
\int_{I_{M}}A_{n}^{\alpha }\left\vert K_{n}^{\alpha }\left( x+t\right)
\right\vert d\mu \left( t\right) \leq \underset{A=0}{\overset{\left\vert
n\right\vert }{\sum }}2^{\alpha A}\int_{I_{M}}\left\vert K_{2^{A}}\left(
x+t\right) \right\vert d\mu \left( t\right) .  \label{9}
\end{equation}

Let $x\in I_{M}\left(e_{k}\right) ,$ $k=0,\dots ,M-1,$ $t\in $ $I_{M}$ \ and 
$x_{q}\neq t_{q},$ where $M\leq q\leq \left\vert n\right\vert -1.$ Using (%
\ref{5a}) we get

\begin{equation*}
\int_{I_{M}}A_{n}^{\alpha }\left\vert K_{n}^{\alpha }\left( x+t\right)
\right\vert d\mu \left( t\right) \leq c_{\alpha }\underset{A=0}{\overset{q-1}%
{\sum }}2^{\alpha A}\int_{I_{M}}2^{k}d\mu \left( t\right) \leq \frac{%
c_{\alpha }2^{k+\alpha q}}{2^{M}}.
\end{equation*}

Hence

\begin{equation}
\int_{I_{M}}\left\vert K_{n}^{\alpha }\left( x+t\right) \right\vert d\mu
\left( t\right) \leq \frac{c_{\alpha }2^{k+\alpha q}}{n^{\alpha }2^{M}}\leq
c_{\alpha }2^{k-M}.  \label{10}
\end{equation}

Let $x\in I_{M}\left(e_{k}\right) ,$ $k=0,\dots ,M-1,$ $t\in $ $I_{M}$ and $%
x_{M}=t_{M},\dots,x_{\left\vert n\right\vert -1}=t_{\left\vert n\right\vert
-1}.$ Applying (\ref{5a}) we have

\begin{equation}
\int_{I_{M}}\left\vert K_{n}^{\alpha }\left( x+t\right) \right\vert d\mu
\left( t\right) \leq \frac{c_{\alpha }}{n^{\alpha }}\overset{\left\vert
n\right\vert -1}{\underset{A=0}{\sum }}2^{\alpha A}\int_{I_{M}}2^{k}d\mu
\left( t\right) \leq c2^{k-M}.  \label{11}
\end{equation}

Combining (\ref{8}), (\ref{10}) and (\ref{11}) we complete the proof of \
Lemma 3.

\section{Proof of the Theorems}

\textbf{Proof of Theorem 1.} By Lemma 1 and (\ref{4}) the proof of first
part of theorem 1 will be complete, if we show that

\begin{equation*}
\int_{I_{M}}\left\vert \overset{\sim }{\sigma }^{\alpha ,\ast
}F(x)\right\vert ^{1/\left(1+\alpha \right) }d\mu \left(x\right) <\infty ,
\end{equation*}
for every $1/\left(1+\alpha \right) $-atom $a.$ We may assume that $a$ be an
arbitrary $1/\left(1+\alpha \right) $-atom with support$\ I,$ $\mu
\left(I\right) =2^{-M}$ and $I=I_{M}.$ It is easy to see that $\sigma
_{n}^{\alpha }\left(a\right) =0,$ when $n\leq 2^{M}.$ Therefore we can
suppose that $n>2^{M}.$

Let $x\in I_{M}.$ Since $\sigma _{n}^{\alpha }$ is bounded from $L_{\infty }$
to $L_{\infty }$ (the boundedness follows from (\ref{4})) and $\left\Vert
a\right\Vert _{\infty }\leq c2^{M/\left( 1+\alpha \right) }$ we obtain 
\begin{eqnarray*}
&&\left\vert \sigma _{n}^{\alpha }a\left( x\right) \right\vert \leq
\int_{I_{M}}\left\vert a\left( t\right) \right\vert \left\vert K_{n}^{\alpha
}\left( x+t\right) \right\vert d\mu \left( t\right) \leq \left\Vert a\left(
x\right) \right\Vert _{\infty }\int_{I_{M}}\left\vert K_{n}^{\alpha }\left(
x+t\right) \right\vert d\mu \left( t\right)  \\
&\leq &c_{\alpha }2^{M\left( 1+\alpha \right) }\int_{I_{M}}\left\vert
K_{n}^{\alpha }\left( x+t\right) \right\vert d\mu \left( t\right) .
\end{eqnarray*}

Let $x\in I_{l+1}\left(e_{k}+e_{l}\right) ,\,0\leq k<l<M.$ From Lemma 3 we
get 
\begin{equation}
\left\vert \sigma _{n}^{\alpha }a\left(x\right) \right\vert \leq \frac{
c_{\alpha }2^{\left(\alpha l+k\right) }2^{\alpha M}}{n^{\alpha }}.
\label{12}
\end{equation}

Let $x\in I_{M}\left(e_{k}\right) ,\,0\leq k<M.$ From Lemma 3 we have

\begin{equation}
\left\vert \sigma _{n}^{\alpha }a\left(x\right) \right\vert \leq c_{\alpha
}2^{\alpha M+k}.  \label{12a}
\end{equation}

Combining (\ref{2}) and (\ref{12}-\ref{12a}) we obtain 
\begin{eqnarray*}
&&\int_{\overline{I_{M}}}\left\vert \overset{\sim }{\sigma }^{\alpha ,\ast
}F(x)\right\vert ^{1/\left(1+\alpha \right) }d\mu \left(x\right) \\
&=&\overset{M-2}{\underset{k=0}{\sum }}\overset{M-1}{\underset{l=k+1}{\sum }}
\int_{I_{l+1}\left(e_{k}+e_{l}\right) }\underset{n>2^{M}}{\sup }\left\vert 
\frac{\sigma _{n}^{\alpha }a\left(x\right) }{\log ^{1+\alpha }n}\right\vert
^{1/\left(1+\alpha \right) }d\mu \left(x\right) \\
&&+\overset{M-1}{\underset{k=0}{\sum }}\int_{I_{M}\left(e_{k}\right) } 
\underset{n>2^{M}}{\sup }\left\vert \frac{\sigma _{n}^{\alpha }a\left(
x\right) }{\log ^{1+\alpha }n}\right\vert ^{1/\left(1+\alpha \right) }d\mu
\left(x\right) \\
&\leq &\frac{1}{M}\overset{M-2}{\underset{k=0}{\sum }}\overset{M-1}{\underset%
{l=k+1}{\sum }}\int_{I_{l+1}\left(e_{k}+e_{l}\right) }\underset{n>2^{M}}{
\sup }\left\vert \sigma _{n}^{\alpha }a\left(x\right) \right\vert
^{1/\left(1+\alpha \right) }d\mu \left(x\right) \\
&&+\frac{1}{M}\overset{M-1}{\underset{k=0}{\sum }}\int_{I_{M}\left(
e_{k}\right) }\underset{n>2^{M}}{\sup }\left\vert \sigma _{n}^{\alpha
}a\left(x\right) \right\vert ^{1/\left(1+\alpha \right) }d\mu \left( x\right)
\\
&\leq &\frac{c_{\alpha }}{M}\overset{M-2}{\underset{k=0}{\sum }}\overset{M-1}%
{\underset{l=k+1}{\sum }}\frac{1}{2^{l}}\frac{2^{\left(\alpha l+k\right)
/\left(1+\alpha \right) }2^{\alpha M/\left(1+\alpha \right) }}{n^{\alpha
/\left(1+\alpha \right) }}+\frac{c_{\alpha }}{M}\overset{M-1}{\underset{k=0}{%
\sum }}\frac{1}{2^{M}}2^{\alpha M/\left(1+\alpha \right) }2^{k/\left(
1+\alpha \right) } \\
&\leq &\frac{c_{\alpha }2^{\alpha M/\left(1+\alpha \right) }}{Mn^{\alpha
/\left(1+\alpha \right) }}\overset{M-2}{\underset{k=0}{\sum }}\overset{M-1}{ 
\underset{l=k+1}{\sum }}\frac{2^{\left(\alpha k+l\right) /\left(1+\alpha
\right) }}{2^{l}}+\frac{c_{\alpha }}{M}\overset{M-1}{\underset{k=0}{\sum }} 
\frac{2^{k/\left(1+\alpha \right) }}{2^{M/\left(1+\alpha \right) }}\leq
c_{\alpha }<\infty .
\end{eqnarray*}

Now, we prove second part of Theorem 1. Let$\ \left\{ \lambda _{k},\text{ }
k\in \mathbb{N}_{+}\right\} $ be an increasing sequence of positive integers
such that 
\begin{equation*}
\overline{\lim_{k\rightarrow \infty }}\frac{\log ^{1+\alpha }\left(\lambda
_{k}\right) }{\varphi \left(\lambda _{k}\right) }=\infty .
\end{equation*}
It is easy to show that for every $\lambda _{k}$ there exists a positive
integer $\ \left\{ n_{k},\text{ }k\in \mathbb{N}_{+}\right\} \subseteq
\left\{ \lambda _{k},\text{ }k\in \mathbb{N}_{+}\right\} $ such that 
\begin{equation*}
\lim_{k\rightarrow \infty }\frac{n_{k}^{1+\alpha }}{\varphi \left(
2^{2n_{k}+1}\right) }=\infty .
\end{equation*}

Let

\begin{equation*}
f_{n_{k}}=D_{2^{2n_{k}+1}}-D_{2^{2n_{k}}}.
\end{equation*}

It is evident 
\begin{equation*}
\widehat{f}_{n_{k}}\left(i\right) =\left\{ 
\begin{array}{l}
\text{ }1,\text{ if }i=2^{2n_{k}},\dots ,2^{2n_{k}+1}-1, \\ 
\text{ }0,\text{otherwise}.%
\end{array}
\right.
\end{equation*}
Then we can write 
\begin{equation}
S_{i}f_{n_{k}}=\left\{ 
\begin{array}{l}
D_{i}-D_{2^{2n_{k}}},\text{ \ \ if \ }i=2^{2n_{k}},\dots ,2^{2n_{k}+1}-1, \\ 
\text{ }f_{n_{k}},\text{ \ \ \ \ \ \ \ \ \ \ \ \ \ if \ }i\geq 2^{2n_{k}+1},
\\ 
0,\text{ \qquad\ \ \ \ \ \ \ \ \ otherwise.}%
\end{array}
\right.  \label{33}
\end{equation}

From (\ref{1dn}) we get 
\begin{equation}
\left\Vert f_{n_{k}}\right\Vert _{H_{1/\left( 1+\alpha \right) }}=\left\Vert
f_{n_{k}}^{\ast }\right\Vert _{H_{1/\left( 1+\alpha \right) }}=\left\Vert
D_{2^{2n_{k}+1}}-D_{2^{2n_{k}}}\right\Vert _{1/\left( 1+\alpha \right) }\leq
c2^{-2\alpha n_{k}}.  \label{34}
\end{equation}%
Let $q_{n_{k}}^{s}=2^{2n_{k}}+2^{2s},$ $s=0,\dots ,n_{k}-1.$ By (\ref{33})
we can write: 
\begin{eqnarray}
&&\frac{\left\vert \sigma _{q_{n_{k}}^{s}}^{\alpha }f_{n_{k}}\right\vert }{%
\varphi \left( q_{n_{k}}^{s}\right) }=\frac{1}{\varphi \left(
q_{n_{k}}^{s}\right) A_{q_{n_{k}}^{s}}^{\alpha }}\left\vert \overset{%
q_{n_{k}}^{s}}{\underset{j=2^{2n_{k}}+1}{\sum }}A_{q_{n_{k}}^{s}-j}^{\alpha
-1}S_{j}f_{n_{k}}\right\vert  \label{35} \\
&=&\frac{1}{\varphi \left( q_{n_{k}}^{s}\right) A_{q_{n_{k}}^{s}}^{\alpha }}%
\left\vert \overset{q_{n_{k}}^{s}}{\underset{j=2^{2n_{k}}+1}{\sum }}%
A_{q_{n_{k}}^{s}-j}^{\alpha -1}\left( D_{j}-D_{2^{2n_{k}}}\right) \right\vert
\notag \\
&=&\frac{1}{\varphi \left( q_{n_{k}}^{s}\right) A_{q_{n_{k}}^{s}}^{\alpha }}%
\left\vert \overset{2^{2s}}{\underset{j=1}{\sum }}A_{2^{2s}-j}^{\alpha
-1}\left( D_{j+2^{2n_{k}}}-D_{2^{2n_{k}}}\right) \right\vert .  \notag
\end{eqnarray}%
Since

\begin{equation}
D_{j+2^{2n_{k}}}-D_{2^{2n_{k}}}=w_{2^{2n_{k}}}D_{j},\text{ \qquad }
j=1,2,\dots,2^{2n_{k}}-1,  \label{36}
\end{equation}
we obtain

\begin{equation}
\frac{\left\vert \sigma _{q_{n_{k}}^{s}}^{\alpha }f_{n_{k}}\right\vert }{
\varphi \left(q_{n_{k}}^{s}\right) }\geq \frac{1}{\varphi \left(
q_{n_{k}}^{s}\right) A_{q_{n_{k}}^{s}}^{\alpha }}\left\vert \overset{2^{2s}}{
\underset{j=0}{\sum }}A_{2^{2s}-j}^{\alpha -1}Dj\right\vert .  \label{37}
\end{equation}

Let $x\in $ $I_{_{2s}}\backslash I_{_{2s+1}}.$ It is easy to show that

\begin{equation}
\frac{\left\vert \sigma _{q_{n_{k}}^{s}}^{\alpha }f_{n_{k}}\left(x\right)
\right\vert }{\varphi \left(q_{n_{k}}^{s}\right) }\geq \frac{
A_{2^{2s}}^{\alpha -1}}{\varphi \left(q_{n_{k}}^{s}\right)
A_{q_{n_{k}}^{s}}^{\alpha }}\overset{2^{2s}}{\underset{j=0}{\sum }}j\geq 
\frac{c2^{4s}A_{2^{2s}}^{\alpha -1}}{\varphi \left(q_{n_{k}}^{s}\right)
A_{q_{n_{k}}^{s}}^{\alpha }}\geq \frac{c2^{2s\left(1+\alpha \right) }}{
\varphi \left(q_{n_{k}}^{s}\right) 2^{2\alpha n_{k}}}.  \label{38}
\end{equation}

Using (\ref{38}) we have 
\begin{equation*}
\int_{G}\left\vert \overset{\sim }{\sigma }^{\alpha ,\ast }f_{n_{k}}\left(
x\right) \right\vert ^{1/\left(1+\alpha \right) }d\mu \left(x\right) \geq 
\text{ }\overset{n_{k}-1}{\underset{s=1}{\sum }}\int_{^{I_{_{2s}}\backslash
I_{_{2s+1}}}}\left\vert \frac{\sigma _{q_{n_{k}}^{s}}^{\alpha
}f_{n_{k}}\left(x\right) }{\varphi \left(q_{n_{k}}^{s}\right) }\right\vert
^{1/\left(1+\alpha \right) }d\mu \left(x\right)
\end{equation*}
\begin{equation*}
\geq c_{\alpha }\overset{n_{k}-1}{\underset{s=1}{\sum }}\frac{2^{2s}}{\left(
\varphi \left(q_{n_{k}}^{s}\right) 2^{2\alpha n_{k}}\right) ^{1/\left(
1+\alpha \right) }}\frac{1}{2^{2s}}\geq \frac{c_{\alpha }n_{k}}{\left(
\varphi \left(2^{2n_{k}+1}\right) 2^{2\alpha n_{k}}\right) ^{1/\left(
1+\alpha \right) }}.
\end{equation*}
From (\ref{34}) we have 
\begin{equation*}
\frac{\left(\int_{G_{m}}\left\vert \overset{\sim }{\sigma }^{\alpha ,\ast
}f_{n_{k}}\right\vert ^{1/\left(1+\alpha \right) }d\mu \right) ^{1+\alpha } 
}{\left\Vert f_{n_{k}}\right\Vert _{H_{1/\left(1+\alpha \right) }}}\geq 
\frac{c_{\alpha }n_{k}^{1+\alpha }}{\varphi \left(2^{2n_{k}+1}\right) }
\rightarrow \infty ,\text{ when }k\rightarrow \infty .
\end{equation*}

Theorem 1 is proved.

\textbf{Proof of Theorem 2. }Suppose that 
\begin{equation*}
\frac{1}{\log n}\overset{n}{\underset{m=1}{\sum }}\frac{\left\Vert \sigma
_{m}^{\alpha }F\right\Vert _{1/\left(1+\alpha \right) }^{1/\left(1+\alpha
\right) }}{m}\leq c_{\alpha }\left\Vert F\right\Vert _{H_{1/\left(1+\alpha
\right) }}^{1/\left(1+\alpha \right) }.
\end{equation*}

Using (\ref{5.1}) we have 
\begin{equation}
\frac{1}{\log n}\overset{n}{\underset{m=1}{\sum }}\frac{\left\Vert \sigma
_{m}^{\alpha }F\right\Vert _{H_{1/\left(1+\alpha \right) }}^{1/\left(
1+\alpha \right) }}{m}=\frac{1}{\log n}\overset{n}{\underset{m=1}{\sum }} 
\frac{\int_{G}\left\Vert \widetilde{\sigma _{m}^{\alpha }F^{\left(t\right) } 
}\right\Vert _{1/\left(1+\alpha \right) }^{1/\left(1+\alpha \right) }dt}{m}
\label{5.3}
\end{equation}
\begin{eqnarray*}
&=&\frac{1}{\log n}\overset{n}{\underset{m=1}{\sum }}\frac{
\int_{G}\left\Vert \sigma _{m}^{\alpha }\widetilde{F^{\left(t\right) }}
\right\Vert _{1/\left(1+\alpha \right) }^{1/\left(1+\alpha \right) }dt}{m}
\leq \int_{G}\frac{1}{\log n}\overset{n}{\underset{m=1}{\sum }}\frac{
\left\Vert \sigma _{m}^{\alpha }\widetilde{F^{\left(t\right) }}\right\Vert
_{1/\left(1+\alpha \right) }^{1/\left(1+\alpha \right) }}{m}dt \\
&\leq &c_{\alpha }\int_{G}\left\Vert \widetilde{F^{\left(t\right) }}
\right\Vert _{H_{1/\left(1+\alpha \right) }}^{1/\left(1+\alpha \right)
}dt\sim c_{\alpha }\int_{G}\left\Vert F\right\Vert _{H_{1/\left(1+\alpha
\right) }}^{1/\left(1+\alpha \right) }dt=c_{\alpha }\left\Vert F\right\Vert
_{H_{1/\left(1+\alpha \right) }}^{1/\left(1+\alpha \right) }.
\end{eqnarray*}

By Theorem W and (\ref{5.3}) the proof of theorem 2 will be complete, if we
show that

\begin{equation*}
\frac{1}{\log n}\overset{n}{\underset{m=1}{\sum }}\frac{\left\Vert \sigma
_{m}^{\alpha }a\right\Vert _{1/\left(1+\alpha \right) }^{1/\left(1+\alpha
\right) }}{m}\leq c_{\alpha }<\infty,
\end{equation*}
for every $1/\left(1+\alpha \right) $-atom $a.$ We may assume that $a$ be an
arbitrary $1/\left(1+\alpha \right) $-atom with support$\ I,$ $\mu
\left(I\right) =2^{-M}$ and $I=I_{M}.$ It is easy to see that $\sigma
_{n}\left(a\right) =0,$ when $n\leq 2^{M}.$ Therefore we can suppose that $%
n>2^{M}.$

Let $x\in I_{M}.$ Since $\sigma _{n}$ is bounded from $L_{\infty }$ to $%
L_{\infty }$ (the boundedness follows from (\ref{4})) and $\left\Vert
a\right\Vert _{\infty }\leq c2^{M/\left(1+\alpha \right) }$ we obtain 
\begin{equation*}
\int_{I_{M}}\left\vert \sigma _{m}^{\alpha }a\left(x\right) \right\vert
^{1/\left(1+\alpha \right) }d\mu \leq \left\Vert a\left(x\right) \right\Vert
_{\infty }^{1/\left(1+\alpha \right) }/2^{M}\leq c_{\alpha }<\infty .
\end{equation*}
Hence 
\begin{equation*}
\frac{1}{\log n}\overset{n}{\underset{m=1}{\sum }}\frac{\int_{I_{M}}\left
\vert \sigma _{m}^{\alpha }a\left(x\right) \right\vert ^{1/\left(1+\alpha
\right) }d\mu }{m}\leq \frac{c_{\alpha }}{\log n}\overset{n}{\underset{m=1}{
\sum }}\frac{1}{m}\leq c_{\alpha }<\infty .
\end{equation*}

Combining (\ref{2}) and (\ref{12}-\ref{12a}) we obtain

\begin{eqnarray*}
&&\frac{1}{\log n}\overset{n}{\underset{m=2^{M}+1}{\sum }}\frac{\int_{ 
\overline{I_{M}}}\left\vert \sigma _{m}^{\alpha }a\left(x\right) \right\vert
^{1/\left(1+\alpha \right) }d\mu \left(x\right) }{m} \\
&=&\frac{1}{\log n}\overset{M-2}{\underset{k=0}{\sum }}\overset{M-1}{ 
\underset{l=k+1}{\sum }}\frac{\int_{I_{l+1}\left(e_{k}+e_{l}\right)
}\left\vert \sigma _{m}^{\alpha }a\left(x\right) \right\vert ^{1/\left(
1+\alpha \right) }d\mu \left(x\right) }{m} \\
&&+\frac{1}{\log n}\overset{n}{\underset{m=2^{M}+1}{\sum }}\overset{M-1}{ 
\underset{k=0}{\sum }}\frac{\int_{I_{M}\left(e_{k}\right) }\left\vert \sigma
_{m}^{\alpha }a\left(x\right) \right\vert ^{1/\left(1+\alpha \right) }d\mu
\left(x\right) }{m} \\
&\leq &\frac{1}{\log n}\left(\overset{n}{\underset{m=2^{M}+1}{\sum }}\frac{
c_{\alpha }2^{\alpha M/\left(1+\alpha \right) }}{m^{\alpha /\left(1+\alpha
\right) +1}}+\overset{n}{\underset{m=2^{M}+1}{\sum }}\frac{c_{\alpha }}{m}
\right) <c_{\alpha }<\infty .
\end{eqnarray*}
which completes the proof of Theorem 2.

\end{document}